\documentclass[letterpaper, 10 pt, conference]{ieeeconf}  

\IEEEoverridecommandlockouts                              
\usepackage{amsmath,amssymb}
\usepackage[hidelinks]{hyperref}
\usepackage{xcolor}
\usepackage{graphicx}
\usepackage{booktabs}
\usepackage{dsfont}
\usepackage{algorithm}
\usepackage{algpseudocode}
\usepackage[style=ieee,backend=bibtex,url=false,doi=false,isbn=false,date=year,citestyle=numeric-comp,maxbibnames=10,minbibnames=10,maxcitenames=10,mincitenames=10]{biblatex}

\usepackage[normalem]{ulem} 

\usepackage{xcolor}

\bibliography{paper}
\usepackage{url}
\newtheorem{rem}{Remark}

\DeclareMathOperator*{\argmin}{argmin} 

\title{\LARGE \bf
Large problems are not necessarily hard: \\
A case study on distributed NMPC paying off
}

\author{Gösta Stomberg, Maurice Raetsch, Alexander Engelmann and Timm Faulwasser
	\thanks{GS and TF are with the Institute of Control Systems, Hamburg University of Technology, 21079 Hamburg, Germany (e-mail: goesta.stomberg@tu-dortmund.de, timm.faulwasser@ieee.org).}
	\thanks{GS was and MR is with the Institute of Energy Systems, Energy Efficiency and Energy Economics, TU Dortmund University, 44227 Dortmund, Germany (e-mail: maurice.raetsch@tu-dortmund.de).}
	\thanks{AE is with logarithmo GmbH \& Co. KG, 44227 Dortmund, Germany (e-mail: alexander.engelmann@ieee.org).}
    \thanks{This work was partially supported by the Deutsche Forschungsgemeinschaft (DFG, German Research Foundation) - project number 527447339.}
}

\begin{document}
	
	\maketitle
	\thispagestyle{empty}
	\pagestyle{empty}

	\begin{abstract} 
	A key motivation in the development of Distributed Model Predictive Control~(DMPC) is to accelerate centralized Model Predictive Control (MPC) for large-scale systems.
	DMPC has the prospect of scaling well by parallelizing computations among subsystems.
	However, communication delays may deteriorate the performance of decentralized optimization, if excessively many iterations are required per control step.
	Moreover, centralized solvers often exhibit faster asymptotic convergence rates and, by parallelizing costly linear algebra operations, they can also benefit from modern multi-core computing architectures.
    On this canvas,	we study the computational performance of cooperative DMPC for linear and nonlinear systems.
	To this end, we apply a tailored decentralized real-time iteration scheme to frequency control for power systems. 
	DMPC scales well for the considered linear and nonlinear benchmarks, as the iteration number does \textit{not} depend on the number of subsystems.
	Comparisons with multi-threaded centralized solvers demonstrate competitive performance of the proposed decentralized optimization algorithms.
	\end{abstract}
	
	\section{Introduction}
	Distributed Model Predictive Control (DMPC) refers to linear and nonlinear MPC variants adapted to graph structures present in Cyber-Physical Systems~(CPS)~\cite{Scattolini2009}. In this paper, the abbreviation DMPC refers to \textit{both} linear and Nonlinear MPC~(NMPC) schemes unless detailed otherwise. 
	Cooperative DMPC schemes promise strong performance by solving a centralized Optimal Control Problem~(OCP) via iterative decentralized optimization.
	Such schemes typically require multiple communication rounds among neighbors in each control step, but avoid central coordination.
	Prior work on cooperative DMPC has considered theoretical guarantees, numerical algorithms, and experimental validation, see, e.g.,~\cite{Giselsson2014,Bestler2019,Stomberg2022,Burk2021b,Stomberg2024b}.
	Besides privacy and resilience,  numerical performance is often quoted as a motivation for DMPC, because expensive computations are parallelized among subsystems.
	However, the scalability of DMPC, i.e., the numerical performance for CPS composed of many subsystems, is not straight forward to ascertain and raises two pivotal issues:
	
	First, it is not clear how DMPC	performs if more subsystem are added to a CPS.
	If the decentralized optimization algorithm parallelizes all steps among subsystems, the execution time per optimizer iteration should not depend on the number of subsystems~\cite{Huber2022}.
	DMPC thus can be expected to scale well, if the addition of new subsystems to a CPS does not increase the required number of optimizer iterations per control step.
	From the theory point of view, the exponential decay of sensitivity in graph-structured Nonlinear Programs~(NLPs) indeed implies that the addition of new subsystems to a CPS has little effect on subsystems far away in the graph~\cite{Shin2022}.
	This effect has been observed in simulations for dual decomposition and the Alternating Direction Method of Multipliers~(ADMM) for convex Quadratic Programs~(QPs)~\cite{Conte2012a}.
	Therein, stronger coupling or more connections among subsystems as well as their instability, yet not the number of subsystems, were found to slow down ADMM convergence.
	
	Second, a crucial  question regarding scalability is whether there are computational merits of DMPC	over centralized MPC.
	Indeed, the latter can also benefit from a wide range of parallelization techniques~\cite{Abugalieh2019,Bernal2024}.
	For instance, parallel linear algebra solvers offer speedup on multi-core CPUs~\cite{Stellato2020,MKLPardiso24}.
	Elaborate numerics parallelize optimizer steps by exploiting application-level sparsity arising for instance in OCPs~\cite{Zavala2008,Word2014}.
	Moreover, implementations for GPUs and FPGAs offer further acceleration~\cite{Pacaud2024,Shin2024,Schubiger2020,Yu2017,Wang2023}.
	Thus, from the computational point of view, it is not clear when DMPC is indeed advantageous over its centralized counterpart.
	
	On this canvas, this paper numerically investigates the scalability of DMPC.
	We consider frequency control for power systems, an application for which DMPC schemes date back to~\cite{Venkat2008} and that has served in prior scalability studies for linear DMPC~\cite{AmoAlonso2023}.
	We extend the case study to nonlinear dynamics and apply a stabilizing decentralized Real-Time-Iteration~(dRTI) scheme~\cite{Stomberg2024}.
	The algorithm executes a fixed number of bi-level decentralized SQP~(dSQP) iterations per control step~\cite{Stomberg2022a} and has been validated for formation control on embedded hardware~\cite{Stomberg2023,Stomberg2024b}.
	
	The contributions of this paper are as follows:
	First, we study the required number of optimizer iterations and solve times for the OCP and analyze the resulting closed-loop control performance.
	Second, we consider nonlinear and linear dynamical systems, i.e., non-convex NLPs and convex QPs, noting that dSQP reduces to ADMM in the latter case.	
	Third, we also test state-of-the-art centralized algorithms combined with multi-threaded linear algebra to provide a performance perspective on the obtained DMPC solve times. 
	
	\textit{Notation:} 
	Given a matrix $A$ and integers $i,j$, $[A]_{ij}$ is the component of $A$ at position $(i,j)$.
	The $i$th component of a vector $a \in \mathbb{R}^n$ is denoted by $[a]_i$.
	The concatenation of vectors $a$ and $b$ into a column vector  is $(a,b)$.
	Given an index set $\mathcal{S}$, $(a_i)_{i \in \mathcal{S}} \doteq (a_1,\dots,a_S)$.
	The horizontal concatenation of matrices $A$ and $B$ is $[A,B]$.
	Given a vector $a$, $A = \mathrm{diag}(a)$ is the diagonal matrix where $[A]_{ii} = [a]_i$.
	The set $\mathbb{I}_{[0,N]}$ with integer $N>0$ denotes the integers in the range $[0,N]$.
	For $a \in \mathbb{R}^n$, the vector $a^+ \in \mathbb{R}^m$ with $m \leq n$ includes the strictly positive components of $a$.
	Likewise, the vector $a^-$ includes the strictly negative components of $a$.

	\section{Problem Statement}\label{sec:problem}
	
	We consider a set $\mathcal{S} = \{1,\dots,S\}$ of coupled systems, each governed by the discrete-time nonlinear dynamics
	\begin{equation}\label{eq:disdyn}
	x_i(t+1) = f_i^\delta(x_i(t),u_i(t),x_{\mathcal{N}_i^\mathrm{in}}(t))), \; x_i(0) = {x_{i,0}}.
	\end{equation}
	For all $i \in \mathcal{S}$, $x_i \in \mathbb{X}_i \subseteq \mathbb{R}^{n_{x,i}}$ is the state of subsystem~$i$, $u_i \in \mathbb{U}_i \subseteq \mathbb{R}^{n_{u,i}}$ is its input, $x_{i,0} \in \mathbb{X}_i$ is the initial condition, $\mathbb{X}_i$ and $\mathbb{U}_i$ are closed sets, and $f_i^\delta: \mathbb{R}^{n_{x,i}} \times \mathbb{R}^{n_{u,i}} \times \mathbb{R}^{n_{x,i}^\mathrm{in}} \rightarrow \mathbb{R}^{n_{x,i}}$.
	The dynamics~\eqref{eq:disdyn} couple each subsystem $i \in \mathcal{S}$ to its in-neighbors $\mathcal{N}_i^\mathrm{in} \subseteq \mathcal{S}$ through the neighboring states $x_{\mathcal{N}_i^\mathrm{in}} \doteq (x_j)_{j \in \mathcal{N}_i^\mathrm{in}}$.
	The sampling interval $\delta > 0$ is equal for all subsystems and each subsystem $i \in \mathcal{S}$ can communicate bi-directionally with all in-neighbors $j \in \mathcal{N}_i^\mathrm{in}$.
	 
	We consider cooperative nonlinear DMPC schemes where, at time $t$, the subsystems solve the OCP
	\begin{subequations}\label{ocp}
		\begin{align}
		\min_{\boldsymbol{x},\boldsymbol{u}} \sum_{i \in \mathcal{S}} \Bigg( \hspace*{-1mm} &\sum_{\tau = 0}^{N-1} \hspace*{-1mm}\ell_i\hspace*{-1mm}\left(x_i[\tau],u_i[\tau],x_{\mathcal{N}_i^\mathrm{in}}[\tau]\right) \hspace*{-1mm} + \hspace*{-1mm} V_{\mathrm{f},i} \left(x_i[N] \right) \hspace*{-1mm} \Bigg) \\
		\nonumber\text{subject} & \text{ to for all } i \in \mathcal{S} \\
		 x_i[\tau+1] &= f_i^\delta(x_i[\tau],u_i[\tau],x_{\mathcal{N}_i^\mathrm{in}}[\tau]) \hspace*{3.1mm} \forall \tau \in \mathbb{I}_{[0,N-1]},\\
		 x_i[0] &= x_i(t),\\
		 x_i[\tau]  &\in \mathbb{X}_i \, \forall \tau \in \mathbb{I}_{[0,N]}, \, u_i[\tau] \in \mathbb{U}_i \, \forall \tau \in \mathbb{I}_{[0,N-1]},\\
		 (x_i[\tau],&x_j[\tau]) \in \mathbb{X}_{ij} \hspace*{0.8cm} \forall j \in \mathcal{N}_i^\mathrm{in}, \; \forall \tau \in \mathbb{I}_{[0,\dots,N]}.\label{eq:xcoup}
		\end{align}
	\end{subequations} 
	We denote predicted variables by square brackets to distinguish between open-loop and closed-loop trajectories. 
	Denote the centralized state and input by $x \doteq (x_i)_{i \in \mathcal{S}} \in \mathbb{R}^{n_x}$ and $u \doteq (u_i)_{i \in \mathcal{S}} \in \mathbb{R}^{n_u}$.
	The decision variables in OCP~\eqref{ocp} are the predicted trajectories over the horizon $N>0$, $\boldsymbol{x} \doteq (x[0], \dots, x[N])$ and $\boldsymbol{u} \doteq (u[0],\dots,u[N-1])$.
	The objective consists of individual stage costs $\ell_i : \mathbb{R}^{{n_{x,i}}} \times \mathbb{R}^{n_{u,i}} \times \mathbb{R}^{n_{x,i}^\mathrm{in}} \rightarrow \mathbb{R}$ and terminal penalties $V_{\mathrm{f},i}: \mathbb{R}^{n_{x,i}} \rightarrow \mathbb{R}$ for all $i \in \mathcal{S}$.
	The constraints~\eqref{eq:xcoup} with the closed sets $\mathbb{X}_{ij}$ couple neighboring states and the extension to coupled input constraints is straight forward and is omitted here for brevity.

	To solve OCP online~\eqref{ocp} via decentralized optimization, we introduce state copies of neighboring subsystems, cf.~\cite{Bestler2019,Stomberg2022}.
	That is, we introduce copies $v_{ij} = x_j$ for all $j \in \mathcal{N}_i^\mathrm{in}$ and replace $x_{\mathcal{N}_i^\mathrm{in}}$ in OCP~\eqref{ocp} by the predicted copy trajectory $\boldsymbol{v}_i \doteq (v_i[0],\dots,v_i[N])$ with $v_i \doteq (v_{ij})_{j \in \mathcal{N}_i^\mathrm{in}}$ for all $i \in \mathcal{S}$.
	Let $z_i \doteq (\boldsymbol{x}_i,\boldsymbol{u}_i,\boldsymbol{v}_i) \in \mathbb{R}^{n_{z,i}}$.	
	We rewrite OCP~\eqref{ocp} as the partially separable NLP
	\begin{subequations}\label{nlp}
		\begin{align}
		\min_{z_1 ,\dots,z_S} &\sum_{i \in \mathcal{S}} f_i (z_i)\label{nlp:obj}\\
		\text{subject to }g_i(z_i) &=0, \quad h_i(z_i) \leq 0 \quad \forall i \in \mathcal{S},\label{nlp:eq}\\ 
		\sum_{i \in \mathcal{S}} E_i z_i &= 0 \label{nlp:coup}
		\end{align}
	\end{subequations}
	with three times continuously differentiable functions $f_i: \mathbb{R}^{n_{z,i}} \rightarrow \mathbb{R}$, $g_i: \mathbb{R}^{n_{z,i}} \rightarrow \mathbb{R}^{n_{g,i}}$, and $h_i: \mathbb{R}^{n_{z,i}} \rightarrow \mathbb{R}^{n_{h,i}}$.
	The sparse matrices $E_i \in \mathbb{R}^{n_c} \times \mathbb{R}^{n_{z,i}}$ couple original and copied states between subsystems and we assume that the matrix $E \doteq [E_1,\dots,E_S]$ has full row rank.
	
	\section{Decentralized Real-Time Iterations}
	
	This section recalls the dRTI scheme from~\cite{Stomberg2024} which, similar to centralized RTIs ~\cite{Gros2020}, applies few SQP iterations per control step to enable real-time execution. 
	Specifically, dRTI deploys a bi-level dSQP scheme which combines an inequality-constrained SQP method on the outer level with ADMM on the inner level, guaranteeing local convergence in the presents of non-convex constraints~\cite{Stomberg2022a}. 
 
	On the outer level, dSQP approximates NLP~\eqref{nlp} at an iterate $z^k \hspace*{-0.5mm} \doteq \hspace*{-0.5mm} (z_1^k,\dots,z_S^k) \hspace*{-0.5mm} \in \mathbb{R}^{n_z}$ as the two-block convex QP
		\begin{align} \label{eq:QPadmm} 
		\min_{\substack{z_i \in \mathbb{Z}_i^k \; \forall i \in \mathcal{S}\\ \bar{z}\in \mathbb{E}} }  \sum_{i \in \mathcal{S}}  f_i^{\mathrm{QP},k}(z_i)
		\textrm{ s.t. } z_i - \bar{z}_i= 0 \; | \; \gamma_i \; ~\forall i \in \mathcal{S}.	
		\end{align}
	The auxiliary decision variables $\bar{z} = (\bar{z}_i)_{i \in \mathcal{S}} \in \mathbb{R}^{n_{z}}$ allow to solve QP~\eqref{eq:QPadmm} via ADMM.
	In~\eqref{eq:QPadmm}, the variable $\gamma_i \in \mathbb{R}^{n_{z,i}}$ is the Lagrange multiplier associated to the consensus constraint of subsystem $i$.
	The objective $f_i^\mathrm{QP,k}(z_i) \doteq (z - z_i^k)^\top H_i^k (z_i - z_i^k)/2 + \nabla f_i(z_i^k)^\top (z_i - z_i^k)$ is convex and the Hessian $H_i^k$ is positive semi-definite for all $i \in \mathcal{S}$. Here, we use the Gauss-Newton approximation $H_i^k \doteq \nabla_{z_iz_i}^2 f_i(z_i^k)$ for quadratic cost functions which often arise in stabilizing MPC schemes~\cite{Gros2020}.
	Note that the matrices $H_i^k$ are constant in this case and can be evaluated offline.
	Alternatively, the exact Hessian of the Lagrangian could be used.
	The subsystem constraints are obtained by linearizing the nonlinear constraints~\eqref{nlp:eq},
	\begin{align*}
	\mathbb{Z}_i^k \doteq \left\{ z_i \in \mathbb{R}^{n_{z,i}} \left| \; \begin{aligned} g_i(z_i^k) + \nabla g_i(z_i^k)^{\top} (z_i-z_i^k)&= 0\\
	h_i(z_i^k) + \nabla h_i(z_i^k)^{\top} (z_i-z_i^k) &\leq 0 \end{aligned} \right. \right\}.
	\end{align*}
	QP~\eqref{eq:QPadmm} couples the subsystems through the set $ \mathbb{E} \doteq \{ \bar{z} \in \mathbb{R}^{n_z} | Ez = 0 \}.$
	Define the augmented Lagrangian to~\eqref{eq:QPadmm} as
	\begin{align*}
	L^k_\rho(z,\bar{z},\gamma) &\doteq \sum_{i \in \mathcal{S}} L_{\rho,i}^k(z_i,\bar{z}_i,\gamma_i),
	\end{align*} where $L_{\rho,i}^k \doteq  f_i^{\mathrm{QP},k}(z_i) + \gamma_i^\top (z_i-\bar{z}_i) + \rho \| z_i - \bar{z}_i \|_2^2 / 2$ with the penalty parameter $\rho > 0$. 
	This allows to solve QP~\eqref{eq:QPadmm} via ADMM on the inner level of dSQP as summarized in Algorithm~\ref{alg:admm}.
	Step~\ref{admm-step:1} of ADMM requires each subsystem to solve a small-scale QP and 	Step~\ref{admm-step2} is an averaging step where the averaging matrix $M_\mathrm{avg} = I - E^\top (E E^\top)^{-1} E$ can be computed offline~\cite{Stomberg2022}. 	
	The dual intialization $\gamma^0 \doteq (\gamma_i^0)_{i \in \mathcal{S}}$ satisfies
	\begin{equation}\label{eq:gam}
	M_\mathrm{avg} \gamma^0 = 0
	\end{equation} to facilitate the averaging, cf.~\cite[Ch. 7]{Boyd2011}.	
 
	Algorithm~\ref{alg:d-SQP} summarizes dSQP, where the derivative computation in Step~\ref{dsqp-stp:1} is parallelized among the subsystems.
	In each control step, dRTI applies $k_\mathrm{max}$ SQP iterations and $l_\mathrm{max}$ ADMM iterations per SQP iteration to meet real-time requirements.
	Crucially, we warm-start $z_i^0$ and $\gamma_i^0$ in dSQP with the solutions obtained in the previous MPC step for all~$i \in \mathcal{S}$.
	
	\begin{rem}[Decentralized implementation] Step~\ref{admm-step2} can be decentralized with neighbor-to-neighbor communication for consensus-type NLPs~\cite{Bestler2019,Stomberg2024,Boyd2011}.
	Here, we instead implement the averaging as a centralized matrix-vector product for execution on one machine.
	Our implementation is thus not subject to communication delays that can dominate the execution time of ADMM~\cite{Burk2021,Stomberg2024b}.
	However, the fundamental trends we will observe in the following sections with respect to dSQP scalability also apply to decentralized implementations, if there is sufficient bandwidth to support a constant per-iteration communication delay between neighbors.	   \hfill {\small{$\Box$}}
	\end{rem}

	\begin{rem}[Closed-loop stability~\cite{Stomberg2024}]
		If the control sampling interval $\delta$ is sufficiently small, if OCP~\eqref{ocp} is stabilizing, if the Hessian of the Lagrangian is used to build QP~\eqref{eq:QPadmm}, and if $l_\mathrm{max}$ is sufficiently large, then the closed-loop system is locally stable~\cite[Thm. 2]{Stomberg2024}.  \hfill {\small{$\Box$}}
	\end{rem}

\begin{algorithm}[t]
	\caption{ADMM for solving QP~
		\eqref{eq:QPadmm}~\cite{Boyd2011}}
	\begin{algorithmic}[1]
		\State Initialization: $\bar{z}_i^0,\gamma_i^0 \, \forall i \in \mathcal{S}$ satisfying~\eqref{eq:gam}, $l_\text{max}$ \label{admm-stp:1}
		\For{ $l = 0, 1, \dots, l_\text{max} - 1$}
		\State $\displaystyle z_i^{l+1} = \argmin_{z_i \in \mathbb{Z}_i^k}  L_{\rho,i}^k(z_i,\bar{z}_i^l,\gamma_i^l)$ for all $i \in \mathcal{S}$\label{admm-step:1}   
		\State $\displaystyle \bar{z}^{l+1} = M_\mathrm{avg} z^{l+1}$  \label{admm-step2}
		\State $\gamma_i^{l+1} = \gamma_i^l + \rho (z_i^{l+1}-\bar{z}_i^{l+1})$ for all $i \in \mathcal{S}$\label{admm-step3}
		\EndFor
		\State \textbf{return} $\bar{z}_i^{l_\text{max}},\gamma_i^{l_\text{max}}$ for all $i \in \mathcal{S}$\label{admm:return}	
	\end{algorithmic} \label{alg:admm}	
\end{algorithm}

\begin{algorithm}[t]
	\caption{dSQP for solving NLP~\eqref{nlp}~\cite{Stomberg2022a}}
	\begin{algorithmic}[1]
		\State Initialization: $z_i^0,\gamma_i^0 \, \forall \,i \in \mathcal{S}$ satisfying~\eqref{eq:gam}, $k_\text{max},l_\text{max}$ \label{dsqp-stp:1}
		\For{$k = 0,1,\dots, k_\text{max}-1$ } for all $i \in \mathcal{S}$ \label{dsqp-step:2}
		\State compute $\nabla f_i^k, g_i^k, \nabla g_i^k, h_i^k, \nabla h_i^k$ and build QP~\eqref{eq:QPadmm} \label{stp:3}
		\State initialize Algorithm~\ref{alg:admm} with $z_i^k,\gamma_i^k,l_\text{max}$ and denote \hspace*{0.4cm} the output by $z_i^{k+1},\gamma_i^{k+1}$\label{dsqp-step:4}
		\EndFor\label{euclidendwhile}
		\State \textbf{return} $z_i^{k_\text{max}}, \gamma_i^{k_\text{max}}$ for all $i \in \mathcal{S}$		
	\end{algorithmic} \label{alg:d-SQP}	
\end{algorithm}	 
		
	\section{Power Network Benchmark and Implementation}
	
\begin{figure}
		\includegraphics[width=\columnwidth]{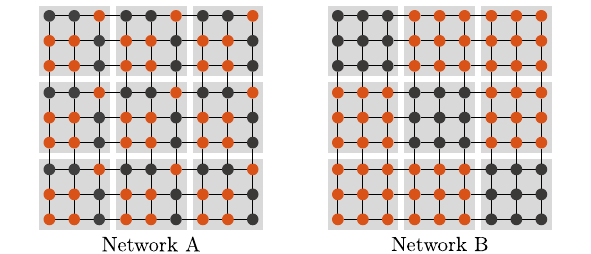}
		\caption{Two 81-bus networks partitioned into nine subsystems of equal size marked in light gray. In Network A, each subsystem contains four generators shown in dark gray and five loads shown in orange. In Network B, three subsystems contain generators and the other subsystems contain loads. Tie lines between coupled buses are shown in black.}\label{fig:grid}
	\end{figure}

	Similar to a prior benchmark on linear DMPC~\cite{AmoAlonso2023}, we consider meshed power networks as shown in Figure~\ref{fig:grid}. 
	The system is well-suited for studying the performance for large-scale OCPs, because more subsystems can be easily added.
	
	\subsection{Power network model}
	
	We consider a bus set $\mathcal{N}$ which we partition into sets of generators $\mathcal{G}$ and loads $\mathcal{L}$.
	For all $n \in \mathcal{N}$, the bus state is $q_n \doteq (\theta_n, \omega_n) \in \mathbb{R}^2$, where the voltage angle $\theta_n$ and angular velocity $\omega_n$ are relative to a synchronous equilibrium, and $\dot{\theta}_n = \omega_n$. 
	The control task is to steer all $\omega_n$ to the origin, i.e., to synchronize all buses at the nominal frequency of $50\,$Hz.	
	The nonlinear synchronous machine dynamics read~\cite{Dorfler2012}	
	\begin{equation}\displaystyle \label{eq:ndyn}
	M_n \dot{\omega}_n(t_{\mathrm{c}})  = -D_n \omega_n(t_{\mathrm{c}}) + P_n(t_{\mathrm{c}}) + p_n(t_{\mathrm{c}}) + w_n(t_{\mathrm{c}})
	\end{equation} for all $n \in \mathcal{N}$, where
	$M_n,D_n>0$ are the inertia and dissipation constants, $t_{\mathrm{c}}$ denotes the continuous time, and
	$P_n(t_{\mathrm{c}}) \doteq - \sum_{m \in \mathcal{M}_n} a_{nm} \sin(\theta_n(t_{\mathrm{c}}) - \theta_m(t_{\mathrm{c}}))$ 
	is the power transfer with neighboring buses.
	The coupling weights are symmetric, $a_{nm} = a_{mn} \geq 0$ for all $n,m \in \mathcal{N}$.
	The set $\mathcal{M}_n \doteq \{m \in \mathcal{N} \; | \; a_{nm} > 0 \}$ collects buses coupled to bus $n \in \mathcal{N}$.
	The controllable power injection at bus $n$ is $p_n \in \mathbb{R}$ and $w_n \in \mathbb{R}$ denotes the uncontrollable load.
	For all generators $n \in \mathcal{G}$, $w_n = 0$ and the generator input constraints are $p_n \in [-0.3,\,0.3]\,$pu.
	Likewise, $p_n = 0$ for all loads $n \in \mathcal{L}$.		
	The state of each bus is constrained to the set~\cite{DeJong2023}
	\begin{equation}\label{eq:ncon}
	\mathbb{Q}_n \hspace*{-1mm}\doteq \hspace*{-1mm}\left\{ q_n \in \mathbb{R}^2 \left | \begin{aligned} &-1.6\pi/\mathrm{s} \leq \omega_n \leq 1.6\pi/\mathrm{s} \\ &-\frac{\pi}{2} \leq \theta_n - \theta_m \leq \frac{\pi}{2} \;  \forall m \in \mathcal{M}_n \end{aligned}\right. \right\}.
	\end{equation} 		
	The nonlinear dynamics~\eqref{eq:ndyn} result in a nonlinear DMPC controller.
	In addition, we consider linear DMPC with linearized dynamics where the power transfer $P_n(t)$ in~\eqref{eq:ndyn} is replaced by
	$
	P_n^\mathrm{lin}(t_{\mathrm{c}}) \doteq - \sum_{m \in \mathcal{M}_n} a_{nm} (\theta_n(t_{\mathrm{c}}) - \theta_m(t_{\mathrm{c}})).
	$
	
	\subsection{Optimal control problem}
	
	To design a DMPC scheme, we partition the network by assigning each bus to one subsystem $i \in \mathcal{S}$.
	Throughout the paper, we consider square networks partitioned into square subsystems of equal size with the meshed coupling structure shown in Figure~\ref{fig:grid}.	
	Let $\mathcal{N}_i$ denote the set of buses assigned to subsystem $i$, let $\mathcal{N}_i \cap \mathcal{N}_j = \emptyset$ for all $i,j \in \mathcal{S}$ if $i \neq j$, and let $\mathcal{N} = \cup_{i \in \mathcal{S}} \mathcal{N}_i$.
	Likewise, we denote the generators and loads assigned to subsystem $i$ by $\mathcal{G}_i$ and $\mathcal{L}_i$, respectively.
	For all $i \in \mathcal{S}$, the subsystem state $x_i \doteq (q_n)_{n \in \mathcal{N}_i}$ includes the states $q_n$ of all buses $n \in \mathcal{N}_i$ and the input $u_i = (p_n)_{n \in \mathcal{G}_i}$ is the controllable power of the generators $n \in \mathcal{G}_i$.
	The uncontrollable loads are modeled as known external disturbances $d_i \doteq (w_n)_{n \in \mathcal{L}_i}$ for all $i \in \mathcal{S}$.
	For all $i \in \mathcal{S}$, the set of in-neighbors is given as 
	\begin{equation*}
	\mathcal{N}_i^\mathrm{in} \doteq \Big\{j \in \mathcal{S} \; \big| \; \exists (n,m) \in \mathcal{N}_i \times \mathcal{N}_j \text{ such that } a_{nm} > 0\Big\}.
	\end{equation*}
	Denote the set of coupled buses from neighboring subsystems as
	$
	\mathcal{M}_i^\mathrm{in} \doteq \cup_{n \in \mathcal{N}_i} \mathcal{M}_n \setminus \mathcal{N}_i.
	$
	For all $i \in \mathcal{S}$, we set the neighboring state to $x_{\mathcal{N}_i^\text{in}} \doteq (\theta_m)_{m \in \mathcal{M}_i^\mathrm{in}}$.
	The careful reader will notice that this choice of $x_{\mathcal{N}_i^\text{in}}$ differs slightly from the definition in Section~\ref{sec:problem}, because we only collect the voltage angles $\theta_m$ of coupled buses $m \in \mathcal{M}_i^\mathrm{in}$ instead of the full neighboring state $x_j$ for all $j \in \mathcal{N}_i^\mathrm{in}$.
	For the considered example, this suffices and incorporating the full neighboring states $x_j$ would unnecessarily inflate the OCP dimension.	
	
	In abstract form, the continuous-time dynamics read
	\begin{equation*}\label{eq:condyn}
	\dot{x}_i(t_{\mathrm{c}}) = f_i^\mathrm{c}(x_i(t_{\mathrm{c}}),u_i(t_{\mathrm{c}}),x_{\mathcal{N}_i^\mathrm{in}}(t_{\mathrm{c}}),d_i(t_{\mathrm{c}})), \; x_i(0) = {x_{i,0}},
	\end{equation*} where, for all $i \in \mathcal{S}$, $f_i^\mathrm{c}: \mathbb{R}^{n_{x,i}} \times \mathbb{R}^{n_{u,i}} \times \mathbb{R}^{n_{x,i}^\mathrm{in}} \times \mathbb{R}^{n_{d,i}} \rightarrow \mathbb{R}^{n_{x,i}}$ includes the nonlinear dynamics~\eqref{eq:ndyn} for all $n \in \mathcal{N}_i$ in the case of nonlinear DMPC or the linear dynamics in the case of linear DMPC.	
	We discretize the continuous-time dynamics with the control sampling interval $\delta$ for all $i \in \mathcal{S}$.
	Specifically, we use zero-order hold discretization for $u_i$ and $d_i$ and we apply the second-order Heun method to those terms in $f_i^\mathrm{c}$ that require only subsystem variables $x_i$, $u_i$, and $d_i$.
	For coupling terms, i.e., $\theta_n - \theta_m$ where $n \in \mathcal{N}_i$ and $m \in \mathcal{M}_i^\mathrm{in}$, we apply the first-order Euler forward discretization when discretizing the neighboring state $\theta_m$. This preserves the sparse coupling structure among subsystems. 
	If a more accurate discretization is desired, tailored distributed multiple shooting techniques could be used~\cite{Savorgnan2011}. 	
	For all $i \in \mathcal{S}$, the state constraint sets $\mathbb{X}_i$ and $\mathbb{X}_{ij}, j \in \mathcal{N}_i^\mathrm{in}$ are composed of the bus constraints $\mathbb{Q}_n$.
	That is, the set $\mathbb{X}_i$ refers to the constraints in~\eqref{eq:ncon} on $\omega_n$ and $\theta_n - \theta_m$ for all $n,m \in \mathcal{N}_i$.
	Likewise, the set $\mathbb{X}_{ij}$ refers to the constraints on $\theta_n - \theta_m$ for all $n \in \mathcal{N}_i$ and $m \in \mathcal{M}_i^\mathrm{in} \cap \mathcal{N}_j$.
	The generator constraints result in a box constraint set $\mathbb{U}_i$ for all $i \in \mathcal{S}$.

	We set $\delta = 0.1\,$s, $N = 100$, $a_{nm} = 0.2\,\text{pu}$ if applicable, and the nominal values $M = 0.167\,\text{pu}\cdot \text{s}^2$ and ${D = 0.045\,\text{pu} \cdot \text{s}}$ similar to~\cite{DeJong2023}.
	For all $i \in \mathcal{S}$, we design quadratic costs
	$
	\ell_i(x_i,u_i) \doteq \sum_{n \in \mathcal{N}_i} \left(q_n^\top Q q_n + R p_n^2 \right)/2
	$ and 
	$ V_{\mathrm{f},i}(x_i) = \sum_{n \in \mathcal{N}_i} q_n^\top P q_n/2 $ with $Q = \mathrm{diag}((0,1))$, $R = 0.1$, and $P = Q$.
	Note that cost functions with coupling between subsystems could also be used if desired.
	Moreover, we regularize the Hessian $H_i$ in the OCP by adding quadratic pentalties with weighting factor $c = 10^{-4}$ to the objective function $f_i$ for all $i \in \mathcal{S}$.

	\subsection{Implementation}

	We have implemented multi-threaded versions of ADMM and dSQP in Julia v. 1.10.5 and the code to perform the simulations is available online.\footnote{\url{https://github.com/OptCon/dmpc_scalability}} 
	All simulations are performed on a Debian Linux Virtual Machine~(VM).
	The VM runs on a server which is equipped with two AMD EPYC 7742 64-core processor sockets and 1$\,$TB RAM.
	The VM has access to 40  CPU cores and 64$\,$GB RAM.
	Furthermore, we test centralized methods equipped with parallel linear algebra solvers to benefit from the multi-core CPU architecture.
	
	As centralized QP solvers for linear-quadratic MPC, we test CPLEX and OSQP v. 0.8.1, referred to as centralized OSQP in the following~\cite{Stellato2020,CPLEX24}.
	For centralized OSQP, we use Intel oneMKL Pardiso as a parallel linear algebra solver~\cite{MKLPardiso24,Schenk2001}.
	To interface CPLEX, we construct and solve a centralized QP model using JuMP v. 1.23.2~\cite{Lubin2023}.
	As decentralized solver, we implement ADMM.
	We use multi threading to accelerate the execution by assigning one thread to each subsystem for Steps~\ref{admm-step:1} and~\ref{admm-step3} of Algorithm~\ref{alg:admm}.
	The averaging matrix $M_\mathrm{avg}$ is computed offline and the subsystem QPs in Step~\ref{admm-step:1} are solved with OSQP running the sequential QDLDL linear algebra solver~\cite{Stellato2020}.
		
	For non-convex NLPs, we use JuMP to model NLP~\eqref{nlp}.
	Specifically, we construct one JuMP model for each subsystem $i \in \mathcal{S}$ to store the objective $f_i$ and the nonlinear constraints $g_i$ and $h_i$.
	This format for storing partially separable NLPs is similar to the Matlab ALADIN-$\alpha$ toolbox~\cite{Engelmann2020}, with the difference that we here use JuMP instead of CasADi for modeling and derivative computation.
	
	We test MadNLP v. 0.8.4 as centralized solver in combination with Intel oneMKL Pardiso~\cite{Shin2020b,Shin2024}.
	As a decentralized method, we test dSQP presented in Algorithm~\ref{alg:d-SQP}, where we evaluate derivatives in Step~\ref{dsqp-stp:1} in parallel across subsystems using automatic differentiation provided by JuMP's MathOptInterface v. 1.32.2~\cite{Legat2021}. 
	For Step~\ref{dsqp-step:4} in dSQP, we use the multi-threaded ADMM implementation for QPs from above.

	\section{Scalability Analysis for Large-Scale OCPs}
	
	This section investigates the optimizer performance on a variety of test cases.  
	Within each case, we increase the network size from 4 to 36 subsystems, but keep the number of decision variables per subsystem constant.

	\subsection{Case study design}
	
	We consider seven test cases grouped into three series as summarized in Table~\ref{tab:parameters}.
	First, we test networks with $9$, $16$, and $25$ buses per subsystem to vary the subsystem size. 
	Second, we change the number of loads per subsystem $|\mathcal{L}_i|$ to adjust the number of control inputs available to each subsystem.
	Third, we consider different initial conditions for the bus frequency $\bar{f}_n \doteq \omega_n / (2\pi)$ for all $n \in \mathcal{N}$.	
	The initial condition for the angles is $\theta_n(0) = 0$ for all $n \in \mathcal{N}$.		

	\begin{table}
		\centering\caption{Problem parameters for the open-loop scalability study.} 
		\begin{tabular}{c c c c l}
			Case & $|\mathcal{N}_i|$ & $|\mathcal{L}_i|$ & $|\bar{f}(0)|$ & Results shown in Figure\\
			\hline
			1 & 9 & 2 & $32\,$mHz & Fig.~\ref{fig:qprw} and Fig.~\ref{fig:nlprw} left columns\\
			2 & 16 & 4 & $32\,$mHz & Fig.~\ref{fig:qprw} and Fig.~\ref{fig:nlprw} center columns\\
			3 & 25 & 5 & $32\,$mHz & Fig.~\ref{fig:qprw} and Fig.~\ref{fig:nlprw} right columns\\
			4 & 9 & 4 & $32\,$mHz & Fig.~\ref{fig:nlpparam} left column\\
			5 & 9 & 6 & $32\,$mHz & Fig.~\ref{fig:nlpparam} left column\\
			6 & 9 & 2 & $48\,$mHz & Fig.~\ref{fig:nlpparam} right column\\
			7 & 9 & 2 & $64\,$mHz & Fig.~\ref{fig:nlpparam} right column\\
		\end{tabular} \label{tab:parameters}
	\end{table}
		
	We proceed as follows for each of the seven test cases.
	First, we choose random parameters for the first subsystem: For all $n \in \mathcal{N}_1$, we sample random parameters from uniform distributions $M_n \in [0.9,1.1]M$, $D_n\in [0.9,1.1]$, and $\bar{f}_n(0) \leq |\bar{f}(0)|$ as given in Table~\ref{tab:parameters}, and we decide at random whether $n \in \mathcal{G}_1$ or $n \in \mathcal{L}_1$.
	Then, we add subsystems with the same parameters to obtain square networks with $4$, $9$, $25$, and $36$ subsystems, see Network A in Figure~\ref{fig:grid}.
	By choosing the same random parameters for all subsystems, we ensure that the only change between different networks within one test case lies in the number of subsystems.
	For OSQP, CPLEX, and MadNLP we report the solve times returned by the solvers.
	For ADMM and dSQP, we first initialize all data structures and then measure the execution times of the code that would be executed in an online DMPC implementation.
	That is, for ADMM and dSQP each subsystem solves the QP in Step~\ref{admm-stp:1} of Algorithm~\ref{alg:admm} once via OSQP before the decentralized solver starts.
	
	For each test case, we tune the penalty parameter $\rho$ in OSQP, ADMM, and dSQP on the network with four subsystems.	
	Moreover, we terminate centralized OSQP, ADMM, and dSQP early based on the centralized KKT residual
	\begin{equation*}
	r \doteq \left\| \left( F_1, \dots, F_s,\sum_{i \in \mathcal{S}} E_i z_i \right) \right \|_\infty, \text{ with} \end{equation*}	
$ F_i \doteq (\nabla_{z_i} L_i, g_i(z_i), h_i^+(z_i),  \mu_i^-, ([\mu_i]_j \cdot [h_i]_j)_{j \in \{1,\dots,n_{h,i}\}})$, multipliers $\nu_i \in \mathbb{R}^{n_{g,i}}, \mu_i \in \mathbb{R}^{n_{h,i}}$, and $\nabla_{z_i} L_i \doteq \nabla f_i(z_i) + \nabla g_i(z_i)\nu_i + \nabla h_i(z_i) \mu_i + \gamma_i$ for all $i \in \mathcal{S}$.

	\subsection{Discussion}
	\subsubsection*{Linear DMPC}
\begin{figure*}
		\includegraphics[width=\textwidth]{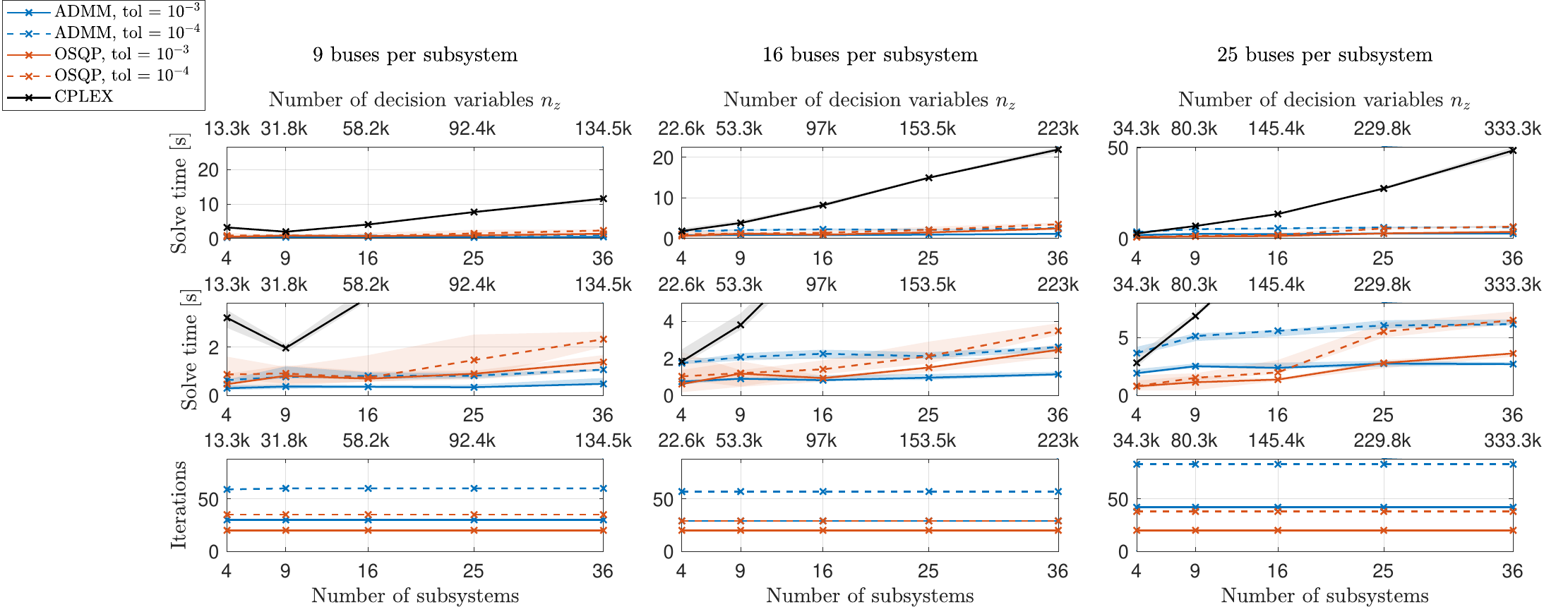}
		\caption{Scalability analysis for convex QPs with $|\mathcal{N}_i| = \{9,16,25\}$ buses per subsystem. The legend on the top left applies to all plots. For each plot, the upper and lower x-axes denote the number of decision variables and the number of subsystems in the network, respectively. Adding more subsystems to the network does not increase the number of necessary iterations, indicating good scalability in decentralized implementations.} \label{fig:qprw}
	\end{figure*}

	Figure~\ref{fig:qprw} shows the performance for convex QPs of cases one to three of Table~\ref{tab:parameters}, including QPs with up to $n_z = 333300$ decision variables.
	The top and center rows display solve times and the center row is a cutout of the top row for better visualization. 
	In Figure~\ref{fig:qprw} and in Figures~\ref{fig:nlprw}--\ref{fig:nlpparam}, crosses mark the median while shaded areas visualize the span from minimum to maximum recorded solve times. 
	CPLEX solves each QP to high accuracy, whereas we terminate ADMM and OSQP early at tolerances of $10^{-3}$ and $10^{-4}$ for the KKT residual.
	Thus, the solve times for CPLEX are not directly comparable to ADMM and OSQP.
	Rather, the CPLEX results show how solving large-scale OCPs to high accuracy becomes computationally infeasible on multi-core CPUs, even with a state-of-the art multi-threaded QP solver.
	In contrast, decentralized ADMM and centralized OSQP rapidly produce suboptimal solutions and scale favorably. 
	Crucially, the required number of centralized OSQP or ADMM iterations does not depend on the number of subsystems, but only on the desired accuracy and the size per subsystem.
	This suggests that both solvers would scale well in practical implementions, as long as the per-iteration time remains constant. 
	We note that for centralized OSQP, tailored implementations for GPUs and FPGAs may provide even further speedup~\cite{Schubiger2020,Wang2023}.
	
	\subsubsection*{Nonlinear DMPC}
\begin{figure*}
		\includegraphics[width=\textwidth]{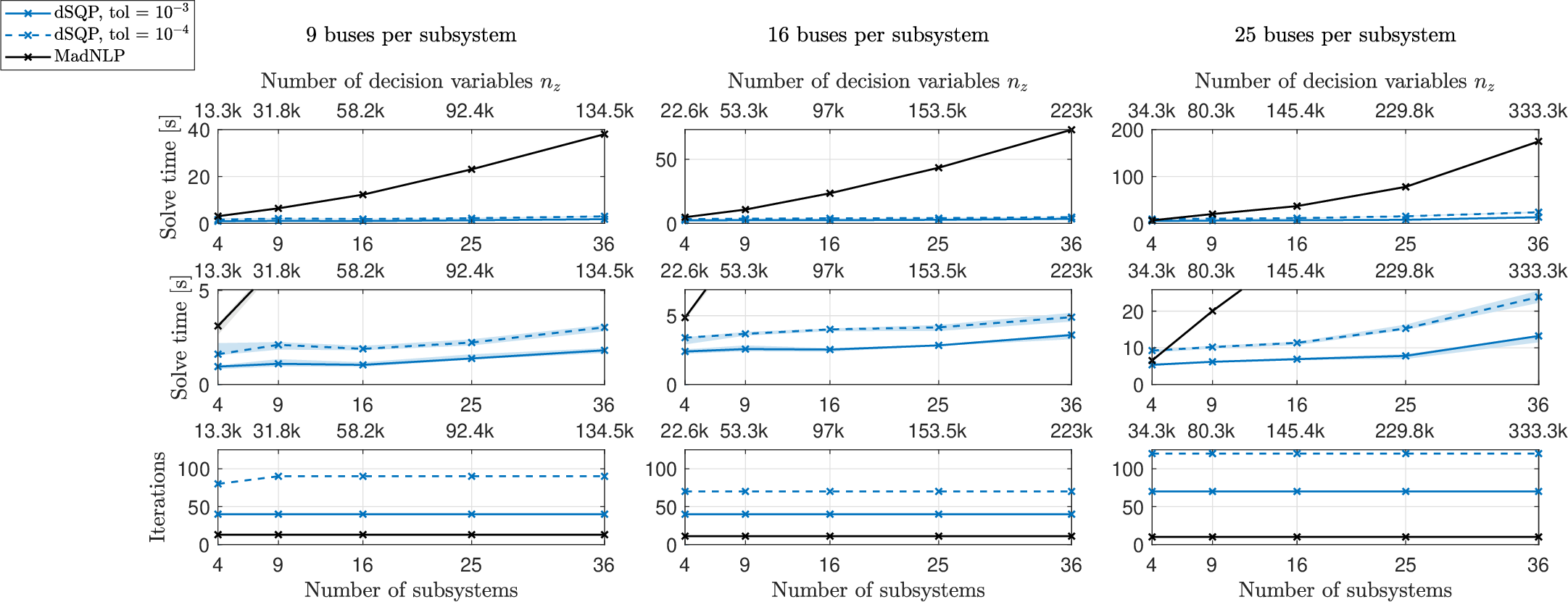}
		\caption{Scalability analysis for non-convex NLPs with varying number of buses per subsystem $|\mathcal{N}_i|$. The legend on the top left applies to all plots. For each plot, the upper and lower x-axes denote the number of decision variables and the number of subsystems in the network, respectively.} \label{fig:nlprw}
	\end{figure*}

    Figure~\ref{fig:nlprw} shows the results for non-convex NLPs with the nonlinear dynamics~\eqref{eq:ndyn} for cases one to three of Table~\ref{tab:parameters}.
	Similarly to CPLEX in the QP simulations, MadNLP solves the NLPs to high accuracy whereas dSQP is terminated early based on the KKT residual.
	The solve times between MadNLP and dSQP are thus not directly comparable as MadNLP achieves greater accuracy. 
	Instead, the MadNLP solve times show the remaining challenges of solving large-scale OCPs to high accuracy with parallelization on CPUs.
	We note that terminating MadNLP early at suboptimal solutions yielded little reduction of computation times. Indeed we observed fast convergence especially in the final interior point iterations and skipping these final iterations had little impact on the overall execution time.
	A promising route to overcome this limitation for centralized MPC is to parallelize centralized solvers on GPUs~\cite{Pacaud2024,Shin2024}.
	
\begin{figure*}
		\centering
		\includegraphics[width=\textwidth]{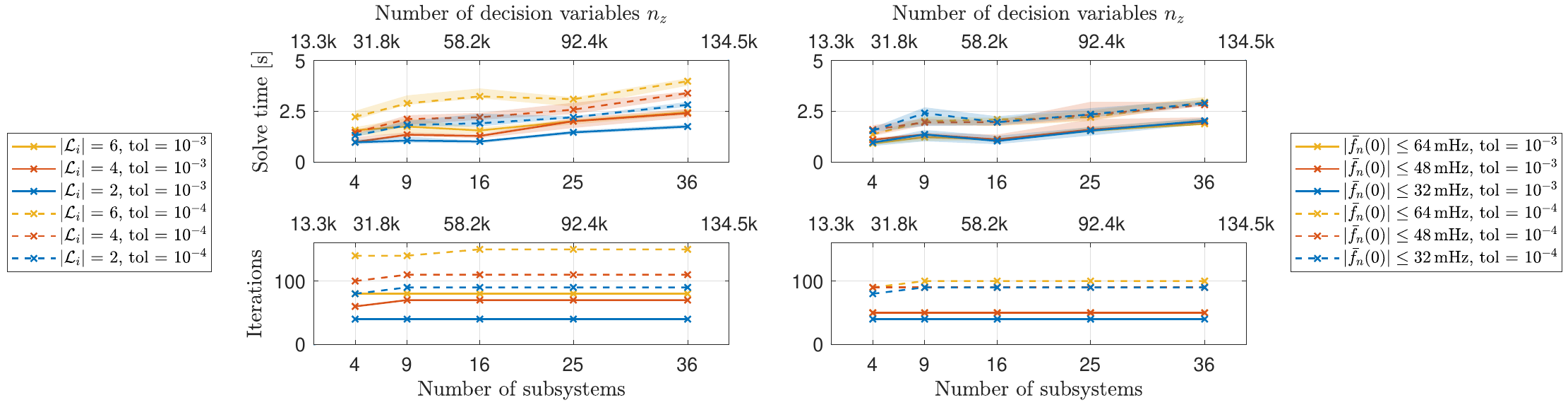}
		\caption{Scalability analysis for dSQP applied to non-convex NLPs with increasing number $| \mathcal{L}_i |$ of loads per subsystem (left) and increasing initial disturbance in the frequency $f_n(0)$ (right). The legend on the left applies to the two plots on the left and the legend on the right applies to the two plots on the right. For each plot, the upper and lower x-axes denote the number of decision variables and the number of subsystems in the network, respectively.} \label{fig:nlpparam}
	\end{figure*}

	Similar to convex QPs, the required number of optimizer iterations shows little dependence on the number of subsystems.
	The same phenomenon can be observed for dSQP in Figure~\ref{fig:nlpparam}, where the number of buses per subsystem is kept constant while the number of loads and the initial condition vary.
	Here, more iterations are required for more challenging subproblems, i.e., with more loads per subsystem, but not with more subsystems in the network.	
	
	\begin{rem}[Distributed vs. centralized MPC?]
	We emphasize that \textit{our} results do not imply that the dRTI scheme outperforms the considered centralized solvers for any given problem.
	For an extensive discussion on problems where ADMM does \textit{not} compare well to centralized solvers, we refer to~\cite{Kozma2015}.
	Instead, our results demonstrate that solving large-scale OCPs to full optimality can be slow.
	We also note that decentralized implementations on multiple machines would see a much slower dRTI execution time due to communication delay, especially for wireless communication~\cite{Stomberg2024b}.
	However, together with the constant per-iteration execution time observed in~\cite{Huber2022}, the constant number of required optimizer iterations for more subsystems suggests good DMPC performance also for large-scale systems.
	In a similar vein, centralized MPC can benefit from this effect.
	We observe constant iteration numbers for OSQP and MadNLP for more subsystems, promising good scalability if the execution time per iteration can be kept low through parallelization.  \hfill {\small{$\Box$}}
	\end{rem}
	
	\section{Closed-loop Control Performance}
	
	ADMM and dSQP scale well for the scenarios in the previos section, because the number of required iterations to reach suboptimal solutions is largely independent of the number of subsystems.
	Hence, this section studies the effect of these suboptimal OCP solutions on closed-loop control performance.
		
	\subsection{Case study design}
	
	We consider the two 81-bus networks shown in Figure~\ref{fig:grid}
	which differ in the positioning of loads and generators.
	In Network A, each subsystem contains both loads and generators and can thus balance power supply and demand. 
	In contrast, only subsystems one, five, and nine of Network~B contain generators whereas the remaining subsystems only contain loads.
	Thus, Network B requires greater consensus among subsystems to balance power.
		
	For both networks, we study a scenario where all buses first start in synchrony and where all loads then exhibit a random perturbation of up to $-0.1\,$pu.	
	To quantify the impact of suboptimal control inputs, we compare the centralized and decentralized solvers from the previous section.
	For linear DMPC and linear centralized MPC, we test $\{0,\dots,10\}$ ADMM or centralized OSQP iterations per control step, respectively.
	For nonlinear DMPC, we use dSQP with $k_\mathrm{max} = 1$ and test the range $l_\mathrm{max} = \{1,\dots,10\}$ of ADMM iterations per SQP iteration.
	Define the averaged closed-loop cost
	\begin{equation*}
	 J \doteq  \frac{1}{t_f} \sum_{t = 0}^{t_n} \sum_{i \in \mathcal{S}} \delta \cdot \ell_i(x_i(t),u_i(t)),
	\end{equation*} where $t_f > 0$ is the duration and $t_n \doteq t_f/\delta$ is the number of MPC steps per simulation.
	We quantify the relative control performance as $J^\star/J$, where $J^\star$ is the averaged closed-loop cost obtained when solving the OCP to optimality in each MPC step via CPLEX or MadNLP, and where $J$ is obtained via centralized OSQP, ADMM, or dSQP.
	
	\begin{figure}
		\includegraphics[width=0.9\columnwidth]{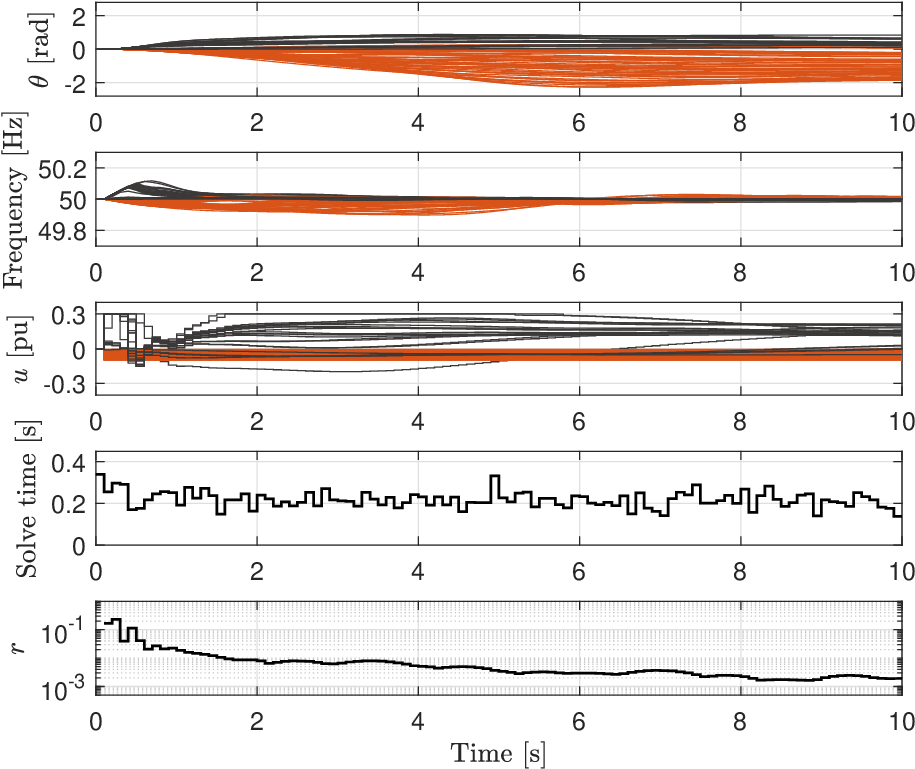}
		\caption{Closed-loop nonlinear DMPC simulation for the 81-bus Network B with $n_z = 31815$ decision variables. The controller runs dSQP with $k_\mathrm{max} = 1$ and $\l_\mathrm{max} = 10$ iterations per control step. Gray lines show generator trajectories and orange lines show load trajectories.}\label{fig:cltraj}
	\end{figure}

	\begin{figure}
		\centering
		\includegraphics[width=0.9\columnwidth]{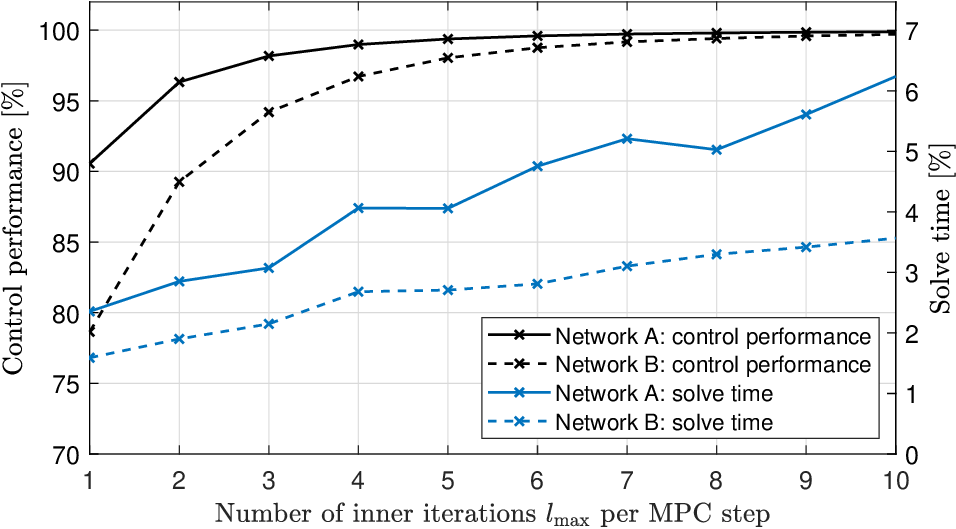}
		\caption{Relative closed-loop control performance $J^\star/J$ for nonlinear DMPC and dSQP solve time as percentage of MadNLP solve time.}\label{fig:clperf}
	\end{figure}
	
	\begin{figure}
		\centering
		\includegraphics[width=0.9\columnwidth]{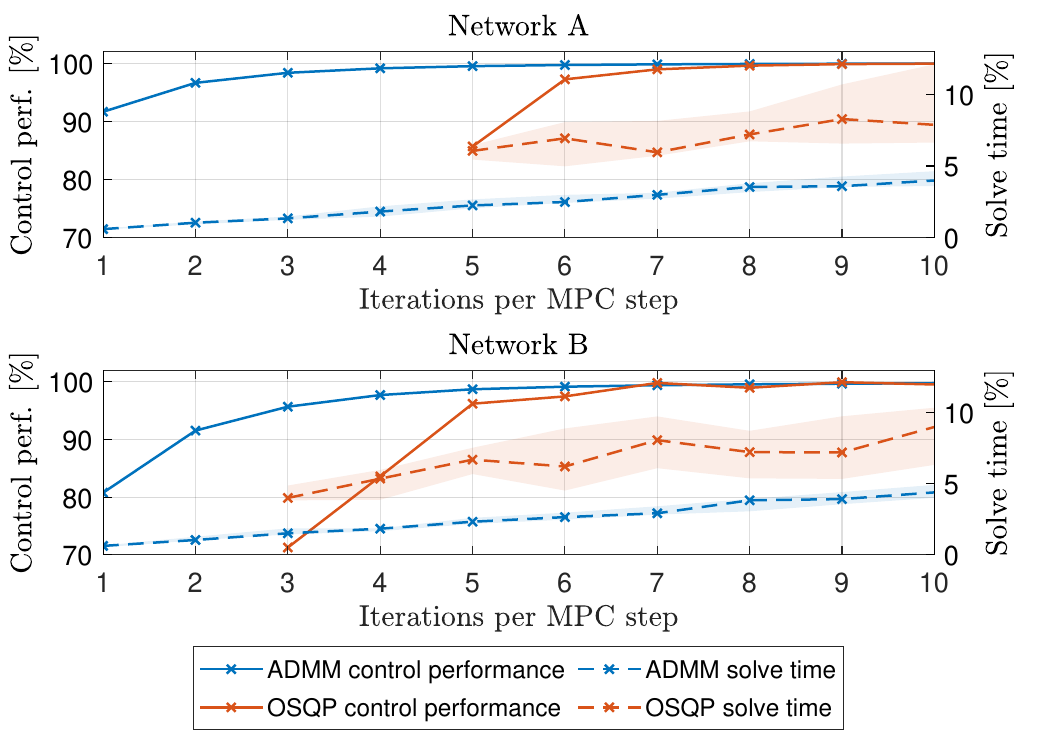}
		\caption{Relative closed-loop control performance $J^\star/J$ for linear DMPC based on ADMM and for linear centralized MPC based on centralized OSQP. Solve times are given as percentage of the CPLEX solve time and the legend on the bottom applies to both plots. For the average solve times per simulation, shaded areas report the minimum and maximum and dashed lines indicate the median.}\label{fig:clperflin}
	\end{figure}
	
	\subsection{Discussion}
	
	\subsubsection*{Nonlinear DMPC}
	Figure~\ref{fig:cltraj} shows closed-loop trajectories for Network $B$ and nonlinear DMPC using dSQP with $l_\mathrm{max} = 10$ iterations per MPC step, less than the 40 iterations required in the previous section for an accuracy of~$10^{-3}$.
	Thus, the accuracy is poor for the first few control steps.
	Nonetheless, dSQP synchronizes the buses thanks to the warm-starting, showcasing the system-optimizer RTI convergence~\cite{Gros2020,Stomberg2024}.
	The center plot in Figure~\ref{fig:cltraj} shows the generator inputs $(p_n)_{n \in \mathcal{G}}$ in dark gray as well as loads in orange, even though $(w_n)_{n \in \mathcal{L}}$ are strictly speaking disturbances and not control signals.
	Observe that the DMPC controls are feasible with respect to the input box constraints, because we only consider state coupling in the OCP. For coupled state constraints, however, the early termination of ADMM can cause constraint violations.		
	Figure~\ref{fig:clperf} shows the relative control performance $J^\star/J$ for $k_\mathrm{max} = 1$ and increasing number of dSQP inner iterations $l_\mathrm{max}$.	
	Furthermore, we plot the dSQP solve time normalized by the MadNLP solve time.
	For both networks, 99$\%$ relative control performance can be obtained with $\l_\mathrm{max} = 7$ ADMM iterations per control step and dSQP requires less than 7\% of the MadNLP execution time.
	The results suggest that DMPC has the potential to scale well as suboptimal inputs in combination with optimizer warm-starting can yield adequate performance.
	At the same time, Network B illustrates a drawback of ADMM-based DMPC: The performance deteriorates for problems that require a high level of consensus unless the number of optimizer iterations per control step is increased.

	\subsubsection*{Linear DMPC}
	Figure~\ref{fig:clperflin} summarizes the control performance and solve times of linear MPC using centralized OSQP and linear DMPC using ADMM.
	The solve times are normalized by the CPLEX solve time.
	The simulations with centralized OSQP and ADMM are run multiple times and Figure~\ref{fig:clperflin} shows the minimum, median, and maximum normalized solve times for each simulation.
	Similar to Figure~\ref{fig:clperf}, suboptimal controls found in a fraction of the high-accuracy solve times suffice for adequate control performance.

	\section{Conclusion}
	This paper has presented simulation results for distributed linear and nonlinear MPC applied to frequency control problems. Decentralized algorithms based on ADMM exhibit good scalability in the considered large-scale examples, because the methods require a constant iteration number to converge as more subsystems are added.
	Future work will examine a wider variety of applications.

	\section*{Acknowledgement}
	The authors thank Emily Elvermann and Tobias Loidl for hardware and software support crucial for the simulations.

	\renewcommand*{\bibfont}{\footnotesize}
	\footnotesize
	\printbibliography

\end{document}